\documentclass{ifacconf}
\usepackage{graphicx} 
\graphicspath{{figures/}} 
\usepackage{natbib}        
\usepackage{amsmath}
\usepackage{amsfonts}
\usepackage{amssymb}
\usepackage{arydshln}
\usepackage{tikz,datatool}
\usepgflibrary{fpu}
\usetikzlibrary{positioning,calc,arrows,shapes,shadows,matrix,backgrounds}
\usepackage{xcolor}

\newcommand{\spa}[1]{\par\vspace*{#1ex}}

\usetikzlibrary{positioning,calc,arrows,shapes,shadows}
\tikzset{
auto,
sys/.style 2 args={
rectangle,
draw,
drop shadow,
fill=white,
minimum height=#2,
minimum width=#1,
inner sep=\dn},
sum/.style={circle,draw,draw=black,inner sep=0mm,minimum size=2mm},
jun/.style={circle,draw,draw=black,inner sep=0mm,minimum size=0mm},
>={latex},
every path/.style={rounded corners},
}

\def\dn{1ex}

\tikzstyle{sy0}=[sys={0*\dn}{0*\dn}]
\tikzstyle{sy1}=[sys={12*\dn}{8*\dn}]
\tikzstyle{sy2}=[sys={8*\dn}{6*\dn}]
\tikzstyle{sy3}=[sys={5*\dn}{5*\dn}]

\newcommand{\sat}{{\rm sat}}
\newcommand{\dzn}{{\rm dzn}}

\newcommand{\Id}{I_d}
\newcommand{\ot}{\otimes}

\newcommand{\arr}[2]{\begin{array}{#1}#2\end{array}}

\newcommand{\te}[1]{\text{\ \ #1\ \ }}

\newcommand{\z}{{\rm z}}

\renewcommand{\r}[1]{(\ref{#1})}
\newcommand{\hl}{\\\hline}

\newcommand{\smat}[2]{{\tiny \left(\begin{array}{#1}#2\end{array}\right)}}
\renewcommand{\c}[1]{{\cal #1}}
\newcommand{\eql}[2]{\begin{equation}\label{#1}#2\end{equation}}
\newcommand{\T}{{\top}}
\newcommand{\Ts}{{\!\top\!}}
\newcommand{\si}{\sigma}
\newcommand{\As}{A}
\newcommand{\Bs}{B}
\newcommand{\Cs}{C}
\newcommand{\Ds}{D}

\newcommand{\Afs}{\c{A}}
\newcommand{\Bfs}{\c{B}}
\newcommand{\Cfs}{\c{C}}
\newcommand{\Dfs}{\c{D}}

\newcommand{\mul}[1]{\begin{multline}#1\end{multline}}
\newcommand{\R}{\mathbb{R}}         
\newcommand{\N}{\mathbb{N}}         
\renewcommand{\S}{\mathbb{S}}          

\newcommand{\I}{\mathbb{I}}

\newcommand{\cl}{\prec}
\newcommand{\cg}{\succ}

\newcommand{\diag}{\mathrm{diag}}
\newcommand{\col}{\mathrm{col}}

\newcommand{\tr}{\mathrm{trace}}

\newcommand{\ga}{\gamma}
\newcommand{\eps}{\varepsilon}
\newcommand{\la}{\lambda}

\newcommand{\al}{\alpha}
\newcommand{\be}{\beta}

\newcommand{\mat}[2]{\left(\begin{array}{@{}#1@{}}#2\end{array}\right)} 

\newenvironment{red_test}{\color{red}}{} 

\newenvironment{blue_test}{\color{blue}}{} 

\renewcommand{\t}{\widetilde}

\newtheorem{theorem}{Theorem}
\newtheorem{lemma}[theorem]{Lemma}
\newtheorem{corollary}[theorem]{Corollary}
\newtheorem{definition}[theorem]{Definition}

\newcommand{\proof}[1]{\begin{proof}#1\end{proof}}

\renewenvironment{proof}[1][Proof]{
	\bf #1. \rm}
{\hfill \footnotesize{$\blacksquare$}\vspace{2ex}}

\renewcommand{\show}[1]{}
\renewcommand{\show}[1]{#1}
\newcommand{\epro}{\mbox{}\hfill\mbox{\rule{2mm}{2mm}} }

\begin{document}
	\begin{frontmatter}
		
		\title{ Robust Exponential Stability and Invariance Guarantees with General Dynamic O'Shea-Zames-Falb Multipliers\thanksref{footnoteinfo} }
		
		\thanks[footnoteinfo]{Funded by Deutsche Forschungsgemeinschaft (DFG, German Research Foundation) under Germany’s Excellence Strategy - EXC 2075 - 390740016. We acknowledge the support by the Stuttgart Center for Simulation Science (SimTech).}
		
		\author[First]{Carsten W. Scherer}
		
		\address[First]{Department of Mathematics, University of Stuttgart, Germany,
			e-mail: carsten.scherer@imng.uni-stuttgart.de}

		\begin{abstract}               
We propose novel time-domain dynamic integral quadratic constraints with a terminal cost for exponentially weighted slope-restricted gradients of not necessarily convex functions. This extends recent results for subdifferentials of convex function and their link to so-called O'Shea-Zames-Falb multipliers. The benefit of
merging time-domain and frequency-domain techniques is demonstrated
for linear saturated systems.
		\end{abstract}

		
	\end{frontmatter}

\section{INTRODUCTION}\label{Sint}

We consider the stability analysis of a loop with  a
discrete-time linear time-invariant (LTI) system
\eql{sys}{x_{t+1}=\As x_t+\Bs w_t,\ z_t=\Cs x_t+\Ds w_t}
(were $A\in\R^{n\times n}$) in feedback with a gradient nonlinearity
\eql{nl}{w_t=\nabla f(z_t)}
with a differentiable function $f:\R^d\to\R$ and for $t\in\N_0$.
It is by now well-established how the stability properties of such interconnections relate to the convergence analysis of optimization algorithms
\citep{LesRec16} or the safety verification of neural networks \citep{FazMor19}.

The general goal in this paper is to characterize robust exponential stability for the state as well as ellipsoidal invariance for the output of \r{sys}
if the gradients are slope-restricted to $[m,L]\subset\R$ \citep{Fre18,GraEbe22,RotGli22a}.
With the left-shift operator $(\sigma x)(t):=x(t+1)$, systems like \r{sys} are
also described as $\si x=Ax+Bw$, $z=Cx+Dw$. Similarly,  \r{nl} reads as
$w=\nabla f(z)$ with $\nabla f(z)_t:=\nabla f(z_t)$ for $t\in\N_0$.

We rely on a discrete-time version of the stability results based on integral quadratic constraints (IQCs) as surveyed in
\cite{Sch22}. This proceeds as follows.
With a filter
\eql{fil0}{
\si \xi=A_\Psi\xi +B_\Psi\mat{c}{z\\w},\ v=C_\Psi\xi+D_\Psi\mat{c}{z\\w},\ \xi_0=0
}
and a symmetric matrix $P$, one assures that the response of \r{fil0} driven by the trajectories of \r{nl} satisfies the IQC
\eql{iqch}{
\sum_{t=0}^{T-1} v_t^\T Pv_t\geq 0\te{for all}T\in\N.
}
In addition, one verifies that the series interconnection of the given system \r{sys} and the filter \r{fil0} is strictly dissipative with respect to the supply rate $v\mapsto -v^\T P v$ and with a positive definite storage function. Classical dissipativity arguments then guarantee the existence of a constant $c$ such that
\eql{stab}{\sum_{t=0}^{T-1} \|x_t\|^2\leq c^2\|x_0\|^2\te{for all}T\in\N,}
uniformly for all trajectories of the interconnection \r{sys}-\r{nl}. If $\Psi$ is the transfer matrix of \r{fil0}, $\Psi^*P\Psi$ is then said to be the dynamic multiplier which assures stability of the loop.

Specifically, we use general
O'Shea-Zames-Falb (OZF) multipliers \citep{WilBro68,FetSch17c,CarHea20}
and their extension to prove exponential stability
\citep{BocLes15,HuSei16,Fre18,MicSch20}.
The paper by \cite{HuSei16} nicely emphasizes  the  discrepancy between time-domain and frequency domain proofs (and the underlying technical delicacies), related to the fact that storage functions as emerging from frequency domain conditions through the application of the Kalman-Yakubovich-Popov Lemma are in general not sign-definite.

The specific technical contribution of this paper is to overcome this discrepancy. We propose novel IQCs for slope-restricted gradients, which extends \r{iqch} to include,  on the right, a nontrivial terminal cost $\xi_T^\T Z\xi_T$ involving the
filter's state at the end time of the interval $[0,T]\cap\N_0$. As argued by \cite{SchVee18} and \cite{Sch22,Sch22a}, this concept leads to a complete (non-conservative) resolution of the dichotomy between
time-domain and frequency-domain IQC results, which involves
a storage function with a positivity condition that is coupled to the terminal cost matrix $Z$.
It is beyond the scope of this paper to prove that our construction is tight.

Another goal of the paper is to build these IQCs in a direct insightful fashion, even for gradients of functions that are not necessarily convex. We avoid technical delicacies about smoothness or invertibility of gradient maps \citep{Fre18}, or about unboundedness of multipliers, which prevents the use of standard IQC results \citep{MicSch20}. This paper extends the results in \citep{Sch22a} that were confined to merely convex functions, and resolves the troubles that emerge in the continuous-time setting as discussed by \cite{FetSch17b}.

One way to investigate exponential convergence properties of feedback loops uses the exponential weighting map
$$
T_\rho(z_0,z_1,z_2,\ldots)=(\rho^0z_0,\rho^1z_1,\rho^2z_2,\ldots).
$$
for some $\rho\in(0,1)$ \citep{DesVid75}. Clearly, $T_\rho$ is linear and invertible with $T_\rho^{-1}=T_{\rho^{-1}}$.
It is easily checked that the set of trajectories $(x,w,z)$ of \r{sys} are in one-to-one correspondence with trajectories
$(\bar x,\bar w,\bar z)$ of 
\eql{sysrho}{
\si\bar x=(\rho^{-1}A)\bar x+(\rho^{-1}B)\bar w,\ \ \bar z=C\bar x+D\bar w
}
under the transformations $\bar x=T_{\rho^{-1}} x$, $\bar w=T_{\rho^{-1}}w$, and $\bar z=T_{\rho^{-1}}z$.
Similarly, \r{nl} translates into
\eql{nlrho}{
\bar w=\Delta^f_\rho(\bar z)}
with the static time-varying operator defined as
\eql{Delrho}{\Delta^f_\rho(\bar z)_t:=(\rho^{-t}\nabla f(\rho^t\bar z_t))_t\te{for}t\in\N_0.}

Ensuring \r{stab} for the transformed loop then guarantees that
the state-trajectory of the original loop \r{sys}-\r{nl} converges exponentially to zero with rate $\rho$.
Additional ellipsoidal invariance properties can be assured by suitable dissipativity arguments.

The remaining paper is structured as follows. In Section~\ref{S1}, we introduce the subdifferential of non-convex slope-restricted functions and a key dissipation inequality. In Section~\ref{S2}, this is used to propose static quadratic constraints for lifted functions, which leads to dynamic IQCs with a non-zero terminal cost as seen in Section~\ref{S3}.
In Section~\ref{S4}, these are used to demonstrate
how to robustly guarantee amplitude bounds for signals in feedback interconnections, while Section~\ref{S5} exhibits the benefit over existing results with a concrete numerical example.

\textbf{Notation.} 
Any tupel $x\in(\R^{n})^p$ is represented as $x=(x_1,\ldots,x_p)$ and $\col(x)=\col(x):=(x_1^\T,\ldots,x_p^\T)^\T$ stacks the entries of the tupel into a column vector. A matrix $A\in\R^{n\times n}$
is doubly hyperdominant (d.h.d.) if all its off-diagonal elements are non-positive
and if $Ae\geq 0$ and $e^\T A\geq 0$ holds, where $e\in\R^n$ is the the all-ones vector \citep{WilBro68}.
For $A,B\in\R^{n\times n}$, $A\cl B$ means that $A$ and $B$ are symmetric and
$B-A$ is positive definite. We denote by $\|x\|^2:=x^\T x$ the Euclidean norm of $x\in\R^d$.
Moreover, $l_{2e}^d$ is the space of all sequences $x:\N_0\to\R^d$, with
$l_2^d$ denoting the subspace of all $x\in l_{2e}^d$ with $\|x\|_2^2:=\sum_{t=0}^\infty \|x_t\|^2<\infty$.
For $h\in\N$, the truncation and lifting operation of some signal $x=(x_0,x_1,\ldots)\in l_{2e}^d$ is defined and denoted as
$x^h:=\col(x_0,\ldots,x_{h-1})$.

\section{On Slope-Restricted Maps}\label{S1}

Recall that $f:\R^d\to\R$ is $m$-strongly convex for $m>0$ if $f-mq$ is convex where
$q(x)=\frac{1}{2}\|x\|^2$. For $d=1$ and if $f$ is two times continuously differentiable, this means that the slope of $f'$ is bounded from below as $f''(x)\geq m$ for all $x\in\R$. It is natural to extend this definition to nonpositive values of $m$ and to upper bounds on the slope as follows.

\definition{\label{Dsec}Let $\si\in\R$ and $f:\R^d\to\R$ be any function. Then $f$ is called $\sigma$-convex if
$f_\sigma:=f-\sigma q$ is convex. Moreover, $f$ is called $\si$-concave if $f^\si:=\si q-f$ is convex.
The set of functions $f:\R^d\to\R$ which are both $m$-convex and $L$-concave is denoted by $\c{S}_{m,L}$.}

For smooth functions, this class has been also considered by \cite{Fre18,GraEbe22,RotGli22a}.
In the sequel, we collect some observations without assumptions on differentiability.

For $\si_1,\si_2\in\R$, we will tacitly use the properties
\eql{help0}{
f_{\si_1}-f_{\si_2}=f_{\si_1}+f^{\si_2}=f^{\si_2}-f^{\si_1}=(\si_2-\si_1)q.
}

If $f\in\c{S}_{m,L}$, then $f_m+f^L=(L-m)q$ is convex, which implies  $L\geq m$.

We extend the notation in Definition~\ref{Dsec} to the sets $\c{S}_{m,\infty}$ (or $\c{S}_{-\infty,L}$) of functions that are merely $m$-convex (or $L$-concave). Then
$\c{S}_{0,\infty}$ just is the set of convex functions and the following properties hold:
$$
\c{S}_{m,\infty}\cap \c{S}_{-\infty,L}=
\c{S}_{m,L}\subset \c{S}_{\si,\mu}\te{if}\si \leq m\leq L\leq \mu.
$$
\show{Indeed, if $f\in \c{S}_{m,L}$ and $\si\leq m$, then
$f_\si$ equals $f_m+(m-\si)q$ and is convex, i.e., $f\in\c{S}_{\si,\infty}$. If $L\leq \mu$, concavity of $f^\mu$ (and thus $f\in\c{S}_{-\infty,\mu}$) follows from $f^\mu=f^L+(\mu-L)q$.}

The definition of the subdifferential $\partial f$ for a convex function $f\in\c{S}_{0,\infty}$ can be seamlessly extended as follows.

\definition{\label{Dsub}For $f\in\c{S}_{m,\infty}$, the subdifferential of $f$ is defined
with any $\si\in(-\infty,m]$ as
\eql{sd1}{\partial f:=
\partial f_\si+\si I.}
For $f\in\c{S}_{-\infty,L}$, it is given for any $\mu\in[L,\infty)$ by
\eql{sd2}{\partial f:=
\mu I-\partial f^\mu        .}
}

On the right in \r{sd1} and \r{sd2}, the usual subdifferential for convex functions appears.
These definitions make sense since they are invariant under the choice of $\si$, $\mu$, respectively.
For \r{sd1} and with $\si_1<\si_2\leq m$, this follows from
$$
\partial f_{\si_1}=\partial (f_{\si_2}+(\si_2-\si_1)q)=\partial f_{\si_2}+(\si_2-\si_1)I,
$$
by using \r{help0} and standard rules for the subdifferential of convex functions. For \r{sd2} and
$L\leq\mu_1<\mu_2$, one relies on
$$
\mu_2I-\partial f^{\mu_2}=\mu_2I-\partial(f^{\mu_1}+(\mu_2-\mu_1)q)=\mu_1I-\partial f^{\mu_1}.
$$

\lemma{For the extended subdifferential, \r{sd1} and \r{sd2} hold as equalities for all $\si\in\R$ and $\mu\in\R$, respectively. Moroever, any function $f\in\c{S}_{m,L}$ with $-\infty<m\leq L<\infty$ is differentiable
and $\partial f$ is equal to the gradient $\nabla f$ of $f$.
}

\begin{proof}
Let $\si>m$. Then $f_\si-(m-\si)q=f_m$ is convex, implying $f_\si\in\c{S}_{m-\si,\infty}$. By \r{sd1} applied to $f_\si$, we infer $\partial f_\si=\partial (f_\si-(m-\si) q)+(m-\si) I=
(\partial f_m+m I)-\si I=\partial f-\si I$. Hence
$\partial f_\si+\si I=\partial f$, which
is \r{sd1} for $\si>m$. The proof of \r{sd2} for $\mu<L$ proceeds analogously.

In case of $f\in\c{S}_{m,L}$ and with the standard subdifferential,  $f_m+f^L=(L-m)q$ implies
$\partial f_m(x)+\partial f^L(x)=(L-m)x$ for $x\in\R^d$. Therefore, $\partial f_m(x)$ is a singleton for every $x\in\R^d$ and, thus, $f_m$ is differentiable. Then the same holds for $f=f_m+mq$, and standard calculus rules
show $\nabla f=\nabla f_\si+\si I=\mu I-\nabla f^\mu$ for all $\si,\mu\in\R$. 
\end{proof}

The following alternative characterizations of vectors in the subdifferentials
in terms of quadratic lower and upper bounding functions extend to any parameters $m,L\in\R$.

\lemma{\label{Lqb}Let $x\in\R^d$. If $f\in\c{S}_{m,\infty}$, then $d\in\partial f(x)$ iff
\eql{sg1}{
f(x)+d^Th+mq(h)\leq f(x+h)\te{for all}h\in\R^d.
}
If $f\in\c{S}_{-\infty,L}$, then $d\in\partial f(x)$ iff
\eql{sg2}{
f(x+h)\leq f(x)+d^Th+Lq(h)\te{for all}h\in\R^d.
}
}
\proof{If $f\in\c{S}_{-\infty,m}$, $f_m$ is convex and by definition of the standard subdifferential, $z\in\partial f_m(x)$ holds iff
\eql{sg}{
f_m(x)+z^Th\leq f_m(x+h)\te{for all}h\in\R^d.
}
Since $q(x+h)-q(x)=x^Th+q(h)$, this is equivalent to
$$
f(x)+(z+mx)^Th+mq(h)\leq f(x+h)\te{for all}h\in\R^d.
$$
The first statement follows from $\partial f_m(x)+mx=\partial f(x)$.

For the second, note that $-f\in \c{S}_{-L,0}\subset\c{S}_{-L,\infty}$.
Hence, $d\in \partial f(x)$ iff $-d\in\partial (-f)(x)$ iff (as just proven)
$$
-f(x+h)\geq -f(x)+(-d)^Th+(-L)q(h)\te{for all}h\in\R^d
$$
iff \r{sg2} is valid.\epro}

From now on we tacitly assume $-\infty <m<L<\infty$
and prove a key inequality for functions in the class $\c{S}_{m,L}$.

\theorem{\label{Tzf}
Let $f\in\c{S}_{m,L}$. Then
\eql{di}{V(u)-V(y)\leq S(u,y)\te{for all}u,y\in\R^d}
is satisfied with
\eql{sf}{V(x):=(L-m)f_m(x)-q(\nabla f_m(x)),}
\eql{sr}{S(u,y):=\nabla f_m(u)^\T[\nabla f^L(u)-\nabla f^L(y)].}
}

\proof{
Let $\al:=L-m>0$. With $u\in\R^d$ and $d:=\nabla f_m(u)$, define the convex function
$$
g(x):=f_m(x)-d^\T x\te{for}x\in\R^n.
$$
Then $g$ is $\al$-concave, since it differs from $f_m$ by an affine function and
$\al q-f_m=(L-m)q-f_m=f^L$ is convex. Hence $g\in\c{S}_{0,\al}$.
For any $y\in\R^d$, Lemma~\ref{Lqb} then implies
$$
\inf_{z\in \R^d}g(z)\leq g(y)+\nabla g(y)^Th+\al q(h)\te{for all}h\in\R^d.
$$
One the one hand, the minimum of the quadratic function in $h$ on the right is
$g(y)-\frac{1}{2\al}\|\nabla g(y)\|^2=g(y)-\frac{1}{\al}q(\nabla g(y))$. On the other hand,
$\nabla g(u)=\nabla f_m(u)-d=0$ implies that the left-hand side equals $g(u)$. Hence
$$
\al g(u)-\al g(y)\leq -q(\nabla g(y)).
$$
With $e:=\nabla f_m(y)$ and since $\nabla g(y)=e-d$, we infer
$$
\al f_m(u)-\al f_m(y)-\al d^T(u-y)\leq -q(e)-d^T(d-e)+q(d).
$$
This gives
$$[\al f_m(u)-q(d)]-[\al f_m(y)-q(e)] \leq d^T[\al(u-y)-(d-e)].$$
The left-hand side just reads $V(u)-V(y)$ by \r{sf}.
The right-hand is $S(u,y)$ in \r{sr}, since $\al(u-y)-(d-e)=
(L-m)(u-y)+m(u-y)-(\nabla f(u)-\nabla f(y))=L(u-y)-(\nabla f(u)-\nabla f(y))$.
This proves \r{di} for \r{sf}-\r{sr}.
\epro
}

We emphasize that \r{di} can be interpreted as a dissipation inequality
as addressed in detail in \cite{Sch22a} for the class $\c{S}_{0,\infty}$.
This is not pursued any further in this paper.

\section{Static Quadratic Constraints}\label{S2}

For stability analysis, we now focus on $f\in\c{S}_{m,L}$ with $\nabla f(0)=0$, the class of which
is denoted by $\c{S}_{m,L}^0$. Note that any $f\in\c{S}_{m,L}^0$ satisfies
$\nabla f_m(0)=0$ and $\nabla f^L(0)=0$.

If $h\in\N$, we further introduce the $h$-lift of $f$ as
\eql{fli}{
\t f(x):=f(x_1)+\cdots+f(x_{h})\te{for}x\in(\R^{d})^h.
}
Since $\nabla \t f(x)=\col(\nabla f(x_1),\ldots,\nabla f(x_{h}))$,
we clearly have $\t f\in\c{S}_{m,L}^0$ and infer that $\nabla \t f$ is diagonally repeated.

Similarly, if defining $\t V$, $\t S$ for $\t f$ as in Theorem~\ref{Tzf}, we get
\eql{dil}{
\t V(u)-\t V(y)\leq \t S(u,y)\te{for all}u,y\in(\R^{d})^h.
}
Note that $\t V$, $\t f_m$, and $\t f^L$ turn out to be the $h$-lifts of $V$, $f_m$, and $f^L$, respectively.

If $P\in\R^{h\times h}$ is any permutation matrix and if using $y:=(P\otimes I_d)u$ in \r{dil},
we can hence conclude $\t V(y)=\t V(u)$ and $\nabla \t f^L(y)=(P\otimes I_d)\nabla  \t f^L(u)$.
Therefore, \r{dil} implies
$$0\leq \nabla \t f_m(u)^\T[(I-P)\otimes I_d]\nabla \t f^L(u)\te{for all}u\in(\R^d)^h.$$
The arguments in \citep{ManSaf05,FetSch17c}  then lead to the following result.

\corollary{\label{Czf}
Let $f\in\c{S}_{m,L}^0$. If $M\in\R^{h\times h}$ is doubly hyperdominant, then
\eql{diP}{0\leq \nabla \t f_m(u)^\Ts (M\otimes I_d)\nabla \t f^L(u)
\te{for all}u\in(\R^d)^h.
}
}

\section{IQCs with Terminal Cost}\label{S3}

On the basis of Corollary~\ref{Czf} we are now in the position to construct so-called dynamic IQCs with a nontrivial terminal cost for the uncertainty $\Delta_\rho^f$ in \r{Delrho}.

First, we construct the filter \r{fil0} driven by $\col(z,w)$. 
With $S_{m,L}:={\arraycolsep.2ex\tiny \mat{cc}{L&-1\\-m&1}}\in\R^{2\times 2}$,
we apply a static transformation
\eql{h00}{
\mat{c}{u_1\\u_2}:=\mat{c}{Lz-w\\-mz+w}=(S_{m,L}\otimes I_d)\mat{c}{z\\w}.
}
This results in
$u_1=\nabla f^L(z)$ and $u_2=\nabla f_m(z)$ in case of
$w=\nabla f(z)$, which motivates \r{h00}.
Then the response $v\!=\!\col(y_1,u_1,y_2,u_2)$ of \r{fil0} is defined by filtering
$u_1$, $u_2$ as
\eql{h01}{
\mat{c}{y_1\\u_1}=\mat{c}{\psi_1\Id\\I_d}u_1\te{and}
\mat{c}{y_2\\u_2}=\mat{c}{\psi_2\Id\\I_d}u_2
}
with the FIR transfer functions
$\psi_1(\z)=\la_0+\sum_{k=1}^{\nu_1}\la_{k}\frac{1}{\z^k}$,  and
$\psi_2(\z)=\sum_{k=1}^{\nu_2}\la_{-k}\frac{1}{\z^k}$, respectively.
The coefficients of these filters of lengths $\nu_1,\nu_2\in\N_0$ are collected as
$$\arraycolsep.5ex
\la:=\mat{ccc|c|ccc}{\la_{\nu_1}&\cdots&\la_{1}&\la_0&\la_{-1}&\ldots&\la_{-\nu_2}}=:\mat{ccc}{\la^1&\la^0&\la^2}.
$$
To construct the resulting overall realization \r{fil0}, let
$J_\nu\in \R^{\nu\times \nu}$
be the upper Jordan block with eigenvalue zero and
$e_1,e_\nu\in\R^\nu$ the first, last standard unit vectors.
Then it is easy to check that the matrices in \r{fil0} can be taken as
$$
A_\Psi:=\mat{cc}{J_{\nu_1}&0\\0&J_{\nu_2}}\ot I_d,\  B_\Psi:=
\left[\mat{cc}{e_{\nu_1}&0\\0&e_{\nu_2}}S_{m,L}\right]\ot I_d,
$$
$$
C_\Psi^\la:=\mat{cc}{\la^1&0\\0&0\hl 0&\la^2\\0&0}\ot I_d,\
D_\Psi^\la:=\left[\mat{cc}{\la^0&0\\1&0\hl 0&0\\0&1}S_{m,L}\right]\ot I_d.
$$

We further introduce the Toeplitz matrix
$T^h(\la)\in\R^{h\times h}$ with the first column
$\col(\la_0,\ldots,\la_{\nu_1},0,\ldots)$
and the first row
$\mat{cccccccccccc}{
\la_0          &\la_{-1}      &\cdots        &\la_{-\nu_2}  &0           &\cdots}$.
For $h=\nu_1+1+\nu_2$ and in the row and column partition $(\nu_1+1)+\nu_2$, $(\nu_2+1)+\nu_1$, respectively, its sub-blocks are denoted as
\eql{Toe2}{
T^{\nu_1+1+\nu_2}(\la)=\mat{cc}{T_{12}(\la)&T_{11}(\la)\\T_{22}(\la)&T_{21}(\la)}.
}
Finally, we set
$F^h_{\rho^{-1}}:=\diag(1,\rho^{-1},\ldots,\rho^{-(h-1)})\in\R^{h\times h}.$

\theorem{\label{Tiqc}Let $f\in\c{S}_{m,L}^0$ for $L>m$ and fix $\rho>0$. Then the response of the filter \r{fil0} driven by $w=\Delta_\rho^f(z)$ satisfies the integral quadratic constraint
\eql{iqc}{
\sum_{t=0}^{T-1} v_t^TPv_t\geq \xi_T^TZ(E)\xi_T\te{for all}T\in\N
}
with the running and terminal cost matrices
$$\arraycolsep0.2ex
P:=\mat{cc}{0&{\arraycolsep.2ex\tiny \mat{cc}{0&1\\1&0}}\!\otimes \! I_{d}\\{\arraycolsep.2ex\tiny \mat{cc}{0&1\\1&0}}\!\otimes\! I_{d}&0},\ Z(E):=\mat{cc}{0&E^\T\!\!\ot\Id\\E\!\ot\!\Id&0},
$$
if the constraints
\eql{lp1}{
F_{\rho^{-1}}^{\nu_1+1+\nu_2}
\mat{cc}{T_{12}(\la)&T_{11}(\la)\\T_{22}(\la)&T_{21}(\la)-E}
F^{\nu_1+1+\nu_2}_{\rho^{-1}}\te{is d.h.d.}
}
and
\eql{lp2}{
\sum_{k=-\nu_2}^{\nu_1}\la_k\rho^k\geq 0,\ \ \sum_{k=-\nu_2}^{\nu_1}\la_k\rho^{-k}\geq 0
}
hold for $\la=(\la^1,\la^0,\la^2)\in\R^{\nu_1+1+\nu_2}$ and $E\in\R^{\nu_2\times \nu_1}$.
}

\proof{
Let $w=\Delta_\rho^f(z)$ for $z\in l_{2e}^d$ and consider the response
$v$ of \r{fil0}. Recall that $v=\col(y_1,u_1,y_2,u_2)$ with the signals
defined by \r{h00}-\r{h01}. We also partition the filter's state-trajectory as $\xi=(\xi_1,\xi_2)$ according to $A_\Psi$.
With the truncation and lifting operation for signals as in the notation section,
it is easy to check that
    \eql{fill}{
\mat    {c}{(\xi_j)_h\\y_j^h\\u_j^h}=
\mat{c}{B_j^h\otimes I_d\\D_j^h\otimes I_d\\I}u_j^h
}
holds, where $B_j^h:=(J_{\nu_j}^{h-1}e_{\nu_j},\ldots,J_{\nu_j}e_{\nu_j},e_{\nu_j})$,
$D_1^h:=T^h(\la^1,\la^0,0)$, and $D_2^h:=T^h(\la^2,0,0)$ (and since
$\xi_0=0$).

With the structure of $P$, the representation \r{fill} shows
\mul{
\frac{1}{2}\sum_{t=0}^{h-1}v_t^\T Pv_t=
\sum_{t=0}^{h-1}\mat{c}{(y_2)_t\\(u_2)_t}^\Ts\! \left({\arraycolsep.2ex \mat{cc}{0&1\\1&0}}\otimes I_d\right)\mat{c}{(y_1)_t\\(u_1)_t}=\\=
\sum_{t=0}^{h-1}\big((y_2)_t^\T\!(u_1)_t\!+\!(u_2)_t^\T\! (y_1)_t\big)=
(y_2^h)^\Ts u_1^h+(u_2^h)^\Ts y_1^h =\\
=(u_2^h)^\T\left[(D_2^h\otimes I_d)^\T+(D_1^h\otimes I_d) \right]u_1^h=\\
=(u_2^h)^\T[T^h(\la)\otimes I_d]u_1^h,
\label{h1}
}
where we exploited $(D_2^h)^\T +D_1^h=T^h(\la)$. Due to the structure of $Z(E)$, we infer
\mul{
\frac{1}{2}\xi_h^\T Z(E)\xi_h=(\xi_2)_h^\T (E\otimes I_d)(\xi_1)_h=\\
=(u_2^h)^\T(B_2^h\otimes I_d)^\T (E\otimes I_d) (B_1^h\otimes I_d)u_1^h=\\
=(u_2^h)^\T\left[((B_2^h)^\T EB_1^h) \otimes I_d \right]u_1^h.
\label{h2}
}
We now exploit $w=\Delta_\rho^f(z)$ to conclude that
\r{h00} leads to $(u_1)_t=\rho^{-t}\nabla f^L(\rho^tz_t)$ and
$(u_2)_t=\rho^{-t}\nabla f_m(\rho^tz_t)$.
With $z_\rho^h:=(F^h_{\rho}\otimes I_d)z^h$, we then get
$$
u_1^h=(F^h_{\rho^{-1}}\otimes I_d)\nabla \t f^L(z_\rho^h),\
u_2^h=(F^h_{\rho^{-1}}\otimes I_d)\nabla \t f_m(z_\rho^h).
$$
If combining with \r{h1}-\r{h2}, \r{iqc} is guaranteed in case that
$$
\nabla \t f_m(z_\rho^T)^\T \Big(M_{\rho^{-1}}^T\otimes I_d\Big)
\nabla \t f^L(z_\rho^T)\geq 0\te{for all}T\in\N,
$$
where we use the abbreviation
\eql{M}{
M_{\rho^{-1}}^h:=F^h_{\rho^{-1}}\Big[T^h(\la)-(B_2^h)^\T EB_1^h\Big]F^h_{\rho^{-1}}.
}
By Corollary~\ref{Czf}, it suffices to show that \r{M} is a d.h.d. matrix in order to conclude the proof.

Let $h=h_0:=\nu_1+1+\nu_2$. By $B_{j}^{h_0}=\mat{cc}{0&I_{\nu_j}}$, observe that
\eql{h02}{
T^{h_0}(\la)-(B_2^{h_0})^\T EB_1^{h_0}=\mat{cc}{T_{12}(\la)&T_{11}(\la)\\T_{22}(\la)&T_{21}(\la)-E}.
}
Hence $M_{\rho^{-1}}^{h_0}$ just equals
\r{lp1} and is, by assumption, d.h.d.

As a first consequence, \r{h02} has non-positive off-diagonal entries,
because the diagonal entries of $F_{\rho^{-1}}^{h_0}$ are positive.
In particular, if recalling \r{Toe2},
we can conclude $\la_k\leq 0$ for $k=-\nu_2,\ldots,-1,1,\ldots,\nu_1$.
With \r{lp2}, we hence infer that
\eql{dhd}{F_{\rho^{-1}}^hT^h(\la)F_{\rho^{-1}}^h\te{is d.h.d. for all}h\in\N.}

Now let $h<h_0$. Then $T^h(\la)-(B_2^h)^\T EB_1^h$ is the right lower $h\times h$ sub-matrix
of $T^{h_0}(\la)-(B_2^{h_0})^\T EB_1^{h_0}$.
Moreover, $F_{\rho^{-1}}^{h}$ is the right-lower $h\times h$ sub-block of $\rho^{h_0-h}F_{\rho^{-1}}^{h_0}$.
Since principal sub-matrices and positive multiples of d.h.d. matrices are d.h.d., we infer that $M_{\rho^{-1}}^{h}$ is a d.h.d. matrix.

Finally let $h>h_0$.  With $\nu:=h-h_0$ and by inspection, we extract the structure
\eql{Th}{\arraycolsep.5ex
T^h(\la)\!-\!(B_2^h)^\T EB_1^h\!=\!
\mat{ccc}{T^\nu(\la)&\leq  0&0\\\leq 0&T_{12}(\la)&T_{11}(\la)\\0&T_{21}(\la)&T_{22}(\la)\!-\!E}
}
in the partition $(\nu+(\nu_1+1)+\nu_2)\times(\nu+(\nu_2+1)+\nu_1)$.
In particular, all off-diagonal entries of \r{Th} are nonpositive, which implies the same for $M_{\rho^{-1}}^h$.
Since $E$ does not affect the first $\nu+\nu_1+1$ rows of \r{Th},
the first $\nu+\nu_1+1$ entries of the column vectors $M^h_{\rho^{-1}}e$ and $F_{\rho^{-1}}^hT^h(\la)F_{\rho^{-1}}^he$ are identical. Due to \r{dhd}, these are nonnegative.
Moreover, by \r{h02} and \r{Th}, each of the last $\nu_2$ entry of $M^h_{\rho^{-1}}e$ differs from
the corresponding one of $M_{\rho^{-1}}^{h_0}e$ by a positive multiple, which shows that it is nonnegative. In total,
we infer $M^h_{\rho^{-1}}e\geq 0$. Since analogous arguments lead to $e^TM^h_{\rho^{-1}}\geq 0$, we have proven that $M^h_{\rho^{-1}}$ is a d.h.d. matrix.
\epro
}

\section{Guaranteeing Stability and Performance}\label{S4}

Now we showcase a typical result that can be formulated with Theorem~\ref{Tiqc} in order to guarantee robust exponential stability and performance for the trajectories of the feedback loop \r{sys}-\r{nl}; from now on we assume $D=0$ which assures that this loop is well-posed.

To this end, we transform \r{sys}-\r{nl} into \r{sysrho}-\r{nlrho} by exponential
signal weighting as discussed in Section~\ref{Sint}. If
the filter \r{fil0} is driven by the signals $\bar z$ and $\bar w$ in the loop,
we obtain a trajectory of
\eql{file}{
\si\eta=\Afs\eta+\Bfs\bar w,\ \ v=\Cfs \eta+\Dfs\bar w,\ \
\bar z=\Cfs_p\eta+\Dfs_p\bar w
}
with $\eta=\col(\xi,\bar x)$
and the series interconnection matrices
\eql{filsym}{
\mat{cc}{\Afs&\Bfs\\\Cfs^\la&\Dfs^\la\\\Cfs_p&0}:=
\mat{cc}{
\mat{cc}{
A_\Psi&B_{\Psi}\smat{c}{\Cs\\0}\\
0     &\rho^{-1}\As         }&
\mat{cc}{
B_{\Psi}\smat{c}{0\\I}\\\rho^{-1}\Bs}\\[2ex]
\mat{cc}{C_\Psi^\la&D_{\Psi}^\la\smat{c}{\Cs\\0}}&
D_{\Psi}^\la\smat{c}{0\\I}\\\mat{cc}{0&C}&0}.
}
Recall the definition of the fixed matrix $P$ and the function $Z(E)$ of $E$ in Theorem~\ref{Tiqc} and introduce the abbreviation
$$
\c{Z}(E):=\mat{cc}{Z(E)&0\\0&0}\in\R^{(n_\Psi+n)\times(n_\Psi+n)}
$$
for a more compact formulation of the following result.

\theorem{\label{Trs}
For fixed $\al,\be\geq 0$, let $\c{X}=\c{X}^\T$, $E$ and $\la$ satisfy \r{lp1}, \r{lp2} and
\eql{lmip}{
\arr{c}{\arraycolsep.3ex
\mat{cc}{\Afs^\T\!\c{X}\Afs\!-\!\c{X}&\Afs^\T\!\c{X}\Bfs\\\Bfs^\T\!\c{X}\Afs&\Bfs^\T\c{X}\Bfs}\!+\!
\mat{cc}{\Cfs^\la&\Dfs^\la\\\Cfs_p&0}^{\!\!\T}\!\!
\mat{cc}{P&0\\0&\al I}\!
\mat{cc}{\Cfs^\la&\Dfs^\la\\\Cfs_p&0}\!\cl 0,
\\[3ex]
\be\Cfs_p^\T \Cfs_p \cl \c{X}+\c{Z}(E).
}
}
Then, for any $f\in\c{S}_{m,L}^0$ and all initial conditions $x_0\in\R^n$, the state trajectory
of \r{sys}-\r{nl} decays exponentially with rate $\rho$ and the following performance condition
is satisfied:
\eql{per}{
\sum_{t=0}^{T-1}\al\|\rho^{-t}z_t\|^2+\be\|\rho^{-T}z_T\|^2\leq x_0^\T Xx_0\ \ \forall\, T\in\N.
}
}

Due to the structure of $P$, we stress that \r{lp1}, \r{lp2} and \r{lmip} are indeed affine in all
optimization variables, including $\la$.

\proof{Any trajectory of the loop \r{sys}-\r{nl} transforms it into one of \r{sysrho}-\r{nlrho} with
$\bar w=\Delta_\rho^f(\bar z)$, which drives \r{fil0} to generate a trajectory of
\r{file}. By slightly perturbing the first LMI in \r{lmip} and
using $v=\Cfs\eta+\Dfs \bar w$ as well as  $\bar z=\Cfs_p\eta$,
right-multiplying $\col(\eta_t,\bar w_t)$ and left-multiplying the transpose and summation for $t=0,\ldots,T-1$
(which are routine dissipation arguments \citep{SchWei11}), we get
$$
\eta_T^\T\c{X}\eta_T-\eta_0^\T\c{X}\eta_0\!+\!\sum_{t=0}^{T-1} v_t^\T\! Pv_t\!+\!
\sum_{t=0}^{T-1}(\al\|\bar z_t\|^2+\eps\|\bar x_t\|^2)
\leq 0
$$
for all $T\in\N$. Note that the perturbation with $\eps$ is introduced to be able to show the exponential convergence of the state-trajectory, as seen below.
We can then exploit \r{iqc} to bound the third term
from below and infer, using
$\eta_0=\col(0,x_0)$, that
\eql{Th1}{
\sum_{t=0}^{T-1} (\al\|\bar z_t\|^2+\eps\|\bar x_t\|^2)+\eta_T^\T[\c{X}+\c{Z}(E)]\eta_T \leq x_0^\T Xx_0
}
for all $T\in\N$. If left- and right-multiplying the second LMI in \r{lmip} with $\eta_T^\T$ and $\eta_T$
and again using $\bar z_T=\Cfs_p\eta_T$, \r{Th1} implies
$\sum_{t=0}^{T-1} (\al\|\bar z_t\|^2+\eps\|\bar x_t\|^2)+\be \|\bar z_T\|^2\leq x_0^\T Xx_0$ for all $T\in\N$. Since $\bar z_t=\rho^{-t}z_t$ and $\bar x_t=\rho^{-t}x_t$,
this proves the exponential decay of $x$ with rate $\rho$ and \r{per}.
\epro}

\section{An Illustrative Numerical Example}\label{S5}

Let us only partially illustrate the benefit of our main results for a scenario with $d=1$.
Specifically, we reveal how to combine time-domain results using generalized sector conditions
with Theorem~\ref{Trs} involving the constructed OZF multipliers from absolute stability theory
in the frequency-domain. We are aware of such attempts undertaken in
\cite{FanLin08}, but the seamless integration into dissipativity theory is new.


For positive $l,L\in\R$ consider the deadzone nonlinearity
\eql{dz}{
\dzn_{l,L}(x)=\left\{\arr{ccc}{0&\text{for}&-l\leq x\leq l,\\L(x-l)&\text{for}&l\leq x,\\L(x+l)&\text{for}&x\leq -l,}\right.
}
which is the gradient of a function in $\c{S}^0_{0,L}$. Therefore, it is possible to use our results for analyzing the stability/performance properties of \r{sys} in feedback with \r{dz}.

Since this nonlinearity is more specific, additional IQCs can be incorporated to potentially improve such tests. A particularly often studied one is based on  the so-called generalized sector condition, as introduced for the saturation function $\sat_{l,L}(x)=Lx-\dzn_{l,L}(x)$; from the very broad range of references, we only mention the books by \cite{HuTee06,TarGar11} for comprehensive discussions. The key is the observation that
\eql{dzi}{\arraycolsep0.2ex
\mat{c}{z-y\\\dzn_{l,L}(z)}^{\!\!\T}\!\! P_L
\mat{c}{z-y\\\dzn_{l,L}(z)}\geq 0
\text{ with }P_L:=\mat{cc}{0&L\\L&-2}
}
holds for all $z,y\in\R$ with $|y|\leq l$. Indeed,
$\dzn_{l,L}(z)$ is a convex combination of $L(z-y)$ and $0$; hence $\dzn_{l,L}(z)=\delta(z-y)$ holds for some $\delta=\delta(y,z)\in[0,L]$; 
then
\r{dzi} follows since
$\col(z-y,\dzn_{l,L}(z))^\T P_L\col(z-y,\dzn_{l,L}(z))=
(z-y)^2\col(1,\delta)^\T P_L\col(1,\delta)=(z-y)^22\delta(L-\delta)\geq 0$.

With the new output $y=Hx$ for \r{sys} and by guaranteeing the bound $\sup_{t\geq 0}|y_t|\leq l$,
\r{dzi} leads to a valid IQC for the signals $z-y$ and $\dzn_{l,L}(z)$ along the loop trajectories.
This can be conically incorporated into the first LMI in \r{lmip} to obtain the subsequent inequality \r{lmi1}, while the amplitude bound on $y$ is guaranteed by \r{lmi3}.

\corollary{\label{Crs}
For fixed $\al,\be\geq 0$ and positive $l,L>0$, let
$\c{X}=\c{X}^\T$, $E$, $\la$, $H$ and $\mu\geq 0$ satisfy the constraints \r{lp1}, \r{lp2} together with
\mul{\label{lmi1}\arraycolsep0.3ex
\mat{cc}{\Afs^\T\c{X}\Afs-\c{X}&\Afs^\T\!\c{X}\Bfs\\\Bfs^\T\!\c{X}\Afs&\Bfs^\T\c{X}\Bfs}\!+\!
\mat{cc}{\Cfs^\la&\Dfs^\la\\\Cfs_p&0}^{\!\!\T}\!\!
\mat{cc}{P&0\\0&\al I}\!
\mat{cc}{\Cfs^\la&\Dfs^\la\\\Cfs_p&0}\!+
\\
+\mu\mat{cc|c}{0&C\!-\!H&0\\0&0&1}^{\!\!\T}\!\!P_L
\mat{cc|c}{0&C\!-\!H&0\\0&0&1}
\cl 0,
}
\spa{-2}
\eql{lmi2}{
\be\Cfs_p^\T \Cfs_p \cl \c{X}+\c{Z}(E),
}
\spa{-3}
\eql{lmi3}{
\mat{cc}{0&H}^\T\mat{cc}{0&H}\cl l^2(\c{X}+\c{Z}(E)).
}
Then the same conclusions as in Theorem~\ref{Trs} can be drawn for the interconnection of \r{sys} with the deadzone nonlinearity \r{dz}, if the system's initial condition satisfies
\eql{ell}{
x_0\in\c{E}_X:=\{\xi\in\R^n\mid \xi^\T X\xi\leq 1\}.}
}

\proof{We follow the dissipativity proof of Theorem~\ref{Trs}. Given any exponentially weighted trajectory of the loop, we first note that \r{lmi3} implies
\eql{Ch1}{
|H\bar x_t|^2\leq l^2\eta_t^\T(\c{X}+\c{Z}(E))\eta_t\te{for all}t\in\N_0.
}

This permits to show $|H\bar x_t|\leq l$ for all $t\in\N_0$ by induction.
Indeed, we have
$\eta_0^\T(\c{X}+\c{Z}(E))\eta_0=
x_0^\T Xx_0\leq 1$ and thus $|H\bar x_0|\leq l$.
Then assume that $|H\bar x_t|\leq l$ for all $t\leq T-1$. For these times $t$, we get
$|\rho^tH\bar x_t|\leq l$ (since $\rho\leq 1$) and we also recall $\bar w_t=\rho^{-t}\dzn_{l,L}(\rho^t\bar z_t)$.
Then \r{dzi} implies
\mul{\label{Ch2}
\mat{c}{\eta_t\\\bar w_t}^{\!\T}\!\! \mat{cc|c}{0&C-H&0\\0&0&1}^{\!\T}\!\!
P_L
\mat{cc|c}{0&C-H&0\\0&0&1}\mat{c}{\eta_t\\\bar w_t}=\\=
\rho^{-2t}\mat{c}{\rho^t\bar z_t-\rho^tH\bar x_t\\\dzn_{l,L}(\rho^t\bar z_t)}^{\!\T}\!\!
P_L
\mat{c}{\rho^t\bar z_t-\rho^tH\bar x_t\\\dzn_{l,L}(\rho^t\bar z_t)}\geq 0
}
for all $t\leq T-1$.
Therefore, the dissipation arguments in the proof of Theorem~\ref{Trs} and applied to
\r{lmi1} lead to \r{Th1}. Since $x_0^\T Xx_0\leq 1$, this
shows $\eta_T^\T(\c{X}+\c{Z}(E))\eta_T\leq 1$ and thus $|H\bar x_T|\leq l$ by \r{Ch1}, which finishes the induction step.

We can now draw the conclusion that \r{Ch2} is valid for all $t\in\N_0$.
Despite the new term in the LMI \r{lmi1} if compared to the first one in \r{lmip}, exactly the same dissipation arguments conclude the proof as for Theorem~\ref{Trs}.
\epro}
\begin{figure}\center
\includegraphics[width=8cm,height=2.7cm]{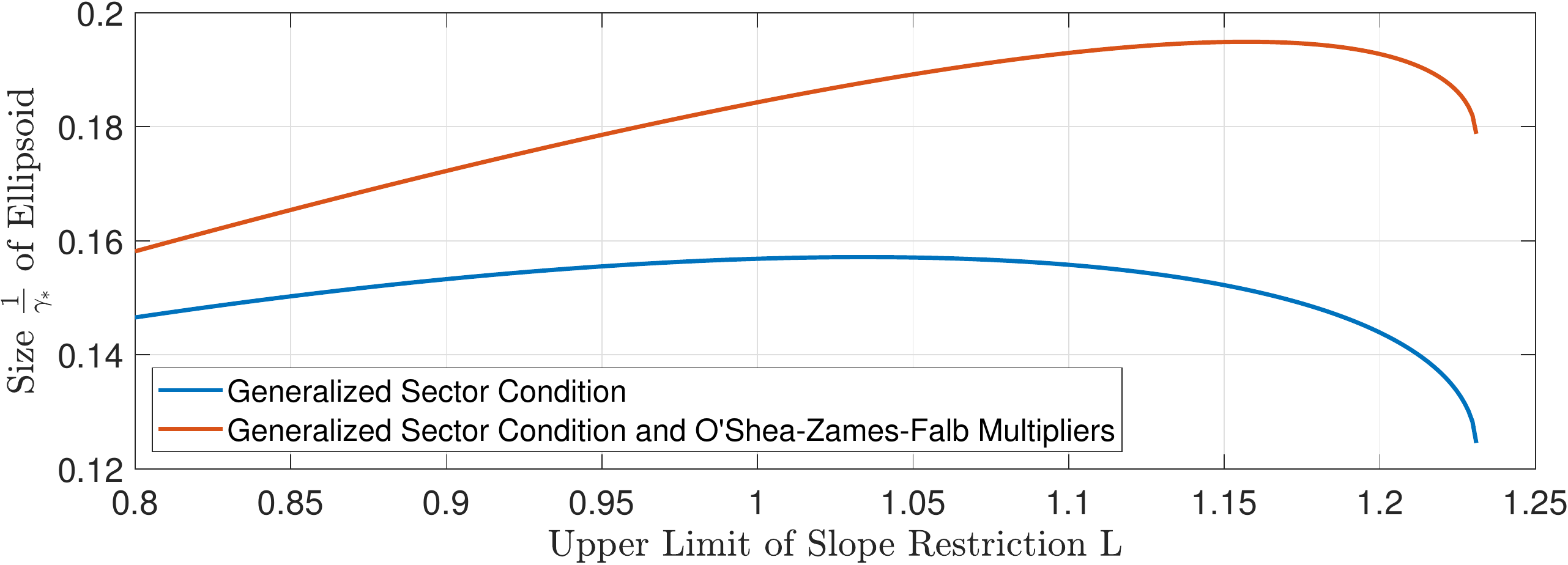}
\caption{Sizes of ellipsoids $\c{E}_{X_*}$ for the generalized sector condition without (blue) and in combination with (red) O'Shea-Zames-Falb multipliers.}\label{fig1}
\end{figure}

For reasons of space, we only exhibit one simple numerical experiment if \r{sys} is defined with
$$
A=\mat{rr}{0.8&0.5\\-0.4&1.2},\ \
B=\mat{c}{-0.18\\1},\ \
C=\mat{rr}{0.3&-1.8}
$$
and $D=0$, interconnected with $w=\sat_{0.1,L}(z)$ for $L\in[0,1.3]$.
We choose $\rho=1$, $\al=0$, $\be=1$ and note  that $1/\sqrt{\mbox{trace}(X)}$ can be considered as a measure for the size of the ellipsoid \r{ell}. This motivates to compute
the infimal $\ga_*$ (with an approximately optimal $X_*\cg 0$) such that the constraints in Corollary \ref{Crs} in addition to $\tr(X)\cl \ga^2I$ hold (by using a line-search over $\mu$ and the LMI-solver of \cite{ML20b} with Yalmip \citep{Lof04}.)

For all $x_0\in \c{E}_{X_*}$, Corollary \ref{Crs} then guarantees the bound
$\sup_{t\geq 0}\|z_t\|\leq 1$.
The sizes $1/\ga_*$ of $\c{E}_{X_*}$ for the generalized sector condition from \citep{HuTee06,TarGar11} depending on $L$ are depicted in
blue in Fig.~\ref{fig1}. The inclusion of dynamic multipliers of length $\nu=\t\nu=1$ leads to an increase of the size of the ellipsoid  (i.e. a reduction of conservatism) as shown by the red curve in Fig.~\ref{fig1}.

\section{Conclusions}\label{Sc}
We have given new time-domain IQCs with a terminal cost for exponentially weighted slope-restricted nonlinearities. For linear saturated systems, it has been demonstrated how these results permit to reduce conservatism by seamlessly merging local time-domain and frequency-domain techniques.
The impact on the analysis of optimization algorithms or the safety verification of linear systems interconnected with neural networks is left for 
future work.

\renewcommand\refname{REFERENCES}

\end{document}